\documentclass[12pt]{amsart}
\usepackage{amssymb,latexsym,cite}
\usepackage[usenames, dvipsnames]{color}
\definecolor{mygray}{gray}{0.6}

\advance\textheight by \topskip   \numberwithin{equation}{section}

\newcommand{\cc}{{\mathbb C}}

\newcommand{\bes}{\begin{split}}
	\newcommand{\fes}{\end{split}}

\newtheorem{theorem}{Theorem}[section]

\newtheorem{corollary}[theorem]{Corollary}

\newtheorem{lemma}[theorem]{Lemma}

\newcommand{\beq}{\begin{eqnarray*}}
\newcommand{\feq}{\end{eqnarray*}}
\newcommand{\beqn}{\begin{eqnarray}}
\newcommand{\feqn}{\end{eqnarray}}

\newlength\myindent
\def\CP{\mathcal{CP}}
\def\CCP{\mathcal{CCP}}
\begin{document}
\title{On column-convex and convex Carlitz polyominoes}

\author{Toufik Mansour}
\address{Department of Mathematics, University of Haifa, 3498838 Haifa, Israel}
\email{tmansour@univ.haifa.ac.il }

\author{Reza Rastegar}
\address{Occidental Petroleum Corporation, Houston, TX 77046 and Departments of Mathematics and Engineering, University of Tulsa, OK 74104, USA - Adjunct Professor}
\email{reza\_rastegar2@oxy.com}

\author{Armend Sh. Shabani}
\address{Department of Mathematics, University of Prishtina, 10000 Prishtin\"{e}, Republic of Kosovo}
\email{armend.shabani@uni-pr.edu }


\thispagestyle{empty}

\begin{abstract}
In this paper, we introduce and study {\it Carlitz polyominoes}. In particular, we show that, as $n$ grows to infinity, asymptotically the number of
\begin{enumerate}
\item column-convex Carlitz polyominoes with perimeter $2n$ is
\beq
\frac{9\sqrt{2}(14+3\sqrt{3})}{2704\sqrt{\pi n^3}}4^n.
\feq
\item convex Carlitz polyominoes with perimeter $2n$ is
\beq
\frac{n+1}{10}\left(\frac{3+\sqrt{5}}{2}\right)^{n-2}.
\feq
\end{enumerate}
\end{abstract}
\maketitle

\section{Introduction and preliminaries}

Take the upright square lattice and call each of its squares, together with its sides, a {\it cell}. A squared lattice polyomino, or simply polyomino, is a finite collection of edge-connected cells in the lattice. The polyominoes were first introduced and studied by Golomb in 1953 \cite{GO} in certain mathematical recreational problems, and soon, they became central objects in the study of Ising model, Pott model, percolation theory, branched polymers \cite{A14, A20, A21, A25}), the mechanics of macromolecules \cite{Tem}, tiling problems \cite{A1, A2, A12, A15, A16, A17}, and enumeration problems \cite{EB, COO, MG}.  We refer to \cite{JG2, A25-2} and references therein for a review of early literature. \par

The enumeration of polyominoes in the general case is an intractable and long-standing open problem. To tackle the problem, many authors have analyzed different subsets of polyominoes assuming additional constraints such as convexity, in which the enumeration has been conducted with respect to several statistics; including but not limited to the area, the perimeter, the outer/inner perimeter, and the corners. The techniques advised in enumeration of these subsets are diverse and range from bijective methods \cite{EB, ADL}, and algorithmic enumeration \cite{JG, COO,JE}, to various forms of decompositions (see \cite{MR} and references wherein). See Table \ref{tab1} for few examples. \par

\begin{table}[htp]
	\begin{tabular}{c||l}
		Description & Ref. \\\hline\hline
		Convex polyominoes & \cite{MG}\\
		Column convex polyominoes &  \cite{FerSvt, Fer1, MR}\\
		$k$-convex polyominoes & \cite{ABS} \\
		Smooth polyominoes & \cite{MSH} \\
		Bargraphs & \cite{Fer1} \\
		Directed animals on a strip & \cite{A14}
	\end{tabular}
	\caption{Various classes of the polyominoes}\label{tab1}
\end{table}

Our objective in this paper is to introduce {\it Carlitz polyominoes} and enumerate them with some additional convexity constraints. To that goal, we first recall some standard definitions. A {\it  column} (resp. {\it row}) of a polyomino is the collection of all cells in the polyomino belonging to an infinite vertical (resp. horizontal) array of  cells in the lattice. For any polyomino $\nu$, let $h(\nu)$ and $v(\nu)$ to be the number of rows and columns of $\nu$, respectively. We define the perimeter of $\nu$ as $2(h(\nu)+v(\nu))$. A polyomino is said to be {\it column-convex} (resp. {\it row convex}) if each of its columns (resp. rows) is a single contiguous block of cells. We say a polyomino is {\it convex} if it is both column-convex and row convex. We further assume that the cells are all squares of size one and equip the lattice with the Cartesian system, in which the bottom-left corner of the bottom cell of the leftmost column at the origin. Numerate the columns from the left side to the right side with the leftmost column counted as ``1". Let $\nu$ be any nonempty polyomino with $m$ columns. We say that the bottom (resp. top) cell of the $i$th column of $\nu$ is at the position $k$ if it lays on (resp. below) and touches the line $y=k$ and we denote this by $b_i(\nu)=k$ (resp. $u_i(\nu)=k$), for $i=1,2,\ldots,m$. Define
\beq
B(\nu)&=&|\{i\mid b_i(\nu)=b_{i+1}(\nu),i=1,2,\ldots,m-1\}|, \\
U(\nu)&=&|\{i\mid u_i(\nu)=u_{i+1}(\nu),i=1,2,\ldots,m-1\}|,
\feq
and refer to $B(\nu)$ and $U(\nu)$ as {\it bottom levels} and {\it top levels}, respectively.

\begin{figure}[htp]
\begin{center}
\begin{picture}(45,28)
\setlength{\unitlength}{.35mm}
\put(30,10){
	\put(20,-10){\line(1,0){20}} \put(10,0){\line(1,0){40}} \put(10,10){\line(1,0){50}} \put(70,10){\line(1,0){10}}
	\put(10,20){\line(1,0){70}} \put(10,30){\line(1,0){80}} \put(10,40){\line(1,0){10}}\put(30,40){\line(1,0){60}}
	\put(40,50){\line(1,0){10}}\put(60,50){\line(1,0){30}}\put(70,60){\line(1,0){20}}
	\put(80,70){\line(1,0){10}}
	\put(10,0){\line(0,1){40}}\put(20,-10){\line(0,1){50}}\put(30,-20){\line(0,1){60}}\put(40,-20){\line(0,1){70}}\put(30,-20){\line(1,0){10}}
	\put(50,0){\line(0,1){50}}\put(60,10){\line(0,1){40}}\put(70,10){\line(0,1){50}}\put(80,10){\line(0,1){60}}\put(90,30){\line(0,1){40}}
}
\end{picture}
\caption{An example of column-convex Carlitz polyomino.}\label{fig 2}
\end{center}
\end{figure}
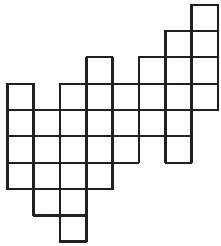

A polyomino $\nu$ is said to be a {\it Carlitz} polyomino if $B(\nu)=U(\nu)=0$.  A polyomino is said to be a convex (resp. column-convex) Carlitz polyomino if it is both convex (resp. column-convex) and Carlitz. See Figure \ref{fig 2} for a column-convex Carlitz polyominoe.

We remark that smooth polyominoes introduced in \cite{MSH} are somewhat related to Carlitz polyominoes. For each smooth polyomino, the top and bottom positions of the columns cannot change by more than one unit going from each column to its neighbors. By contrast, the top and bottom positions in a Carlitz polyomino changes by at least one unit going from each column to its neighbors.

The paper is organized as follows. Section \ref{sec-a} is devoted to the enumeration of the column-convex Carlitz polyominoes. We first obtain an explicit formula for the generating function of the number of these polyominoes according to their perimeter. Then, Corollary \ref{cor2.4} yields that the number of column-convex Carlitz polyominoes with the perimeter $2n$ is asymptotic to
\beq
\frac{9\sqrt{2}(14+3\sqrt{3})}{2704\sqrt{\pi n^3}}4^n,
\feq
as $n$ grows to infinity.
In Section \ref{sec-b}, the convex Carlitz polyominoes are counted. We obtain the generating function for the number of convex Carlitz  polyominoes according to their perimeters. Then, the explicit form is used to show that, as $n$ goes to infinity, the number of convex Carlitz polyominoes of the perimeter $2n$ is asymptotic to
\beq
\frac{n+1}{10}\left(\frac{3+\sqrt{5}}{2}\right)^{n-2}.
\feq

\section{Enumeration of column-convex Carlitz polyominoes} \label{sec-a}

In this section, we consider certain decompositions of the column-convex polyominoes and establish a functional system of linear combinations of multivariate series of forms $A(t_1,\ldots, t_r;x_1,\ldots,x_s)$. These forms do not depend on $x_1,\ldots,x_s,$ simultaneously and, their coefficients are referred to as kernels \cite{doron0}. A systematic approach to solve these type of equations is
given in \cite{Ban}. We conduct all of our calculations and manipulations having the following facts in mind. Let $\mathbb{Q}$, $\mathbb{Q}[x_1, \dots, x_s],$ and $\mathbb{Q}[[x_1, \ldots, x_s]]$ denote, respectively, the field of rational numbers, the ring of polynomials in $x_1,\ldots,x_s$, and the ring of formal power
series in $x_1,\ldots,x_s$ with coefficients in $\mathbb{Q}$.
Recall that a series $A\in \mathbb{Q}[x_1, \dots, x_s][[t_1,\ldots, t_r]]$ is $D$-finite if its partial derivatives span a finite dimensional vector space
over the field of rational functions in $t_1,\ldots,t_r$ with coefficients in
$\mathbb{Q}[x_1,\ldots,x_s]$. Note that any algebraic series is $D$-finite. The specializations of a $D$-finite series, obtained by assigning values in $\mathbb{Q}$ to a subset of variables,
are $D$-finite, if well-defined.  Moreover, if $A$ is $D$-finite, then
any substitution $x_i=1$ with $i\in I\subseteq [s]$ or/and $x_i=x_j$
for $(i,j)\in I\times I\subseteq [s]\times[s]$ into $A$ is also
$D$-finite. One last remark is that to derive asymptotic forms for our results, we apply singularity analysis (see \cite[Section VI]{FS} for a comprehensive review). All calculations for asymptotic analysis require an appropriate domain of complex number $\cc$. We omit the details for the sake of brevity as it is standard.

\begin{figure}[htp]
\begin{picture}(100,26)(0,-10)
\setlength{\unitlength}{.3mm}
\put(0,0){\put(0,0){\line(1,0){10}}\put(10,0){\line(0,1){50}}
\put(10,50){\line(-1,0){10}}\put(0,50){\line(0,-1){50}}
\put(1,-10){\tiny$a$}}
\put(50,0){\put(0,0){\line(1,0){10}}\put(10,0){\line(0,1){50}}
\put(10,50){\line(-1,0){10}}\put(0,50){\line(0,-1){50}}
\put(1,-10){\tiny$a$}
\put(10,10){\put(0,0){\line(1,0){10}}\put(10,0){\line(0,1){30}}
\put(10,30){\line(-1,0){10}}\put(0,30){\line(0,-1){30}}
\put(2,-20){\tiny$s$}\put(18,10){\tiny$\cdots$}}
\put(-15,-35){\tiny$1\leq s\leq a$}}
\put(140,0){\put(0,0){\line(1,0){10}}\put(10,0){\line(0,1){50}}
\put(10,50){\line(-1,0){10}}\put(0,50){\line(0,-1){50}}
\put(1,-10){\tiny$a$}
\put(10,-20){\put(0,0){\line(1,0){10}}\put(10,0){\line(0,1){30}}
\put(10,30){\line(-1,0){10}}\put(0,30){\line(0,-1){30}}
\put(-15,-15){\tiny$s\geq2$}\put(18,10){\tiny$\cdots$}}}
\put(240,0){\put(0,0){\line(1,0){10}}\put(10,0){\line(0,1){50}}
\put(10,50){\line(-1,0){10}}\put(0,50){\line(0,-1){50}}
\put(1,-10){\tiny$a$}
\put(10,30){\put(0,0){\line(1,0){10}}\put(10,0){\line(0,1){30}}
\put(10,30){\line(-1,0){10}}\put(0,30){\line(0,-1){30}}
\put(-15,-65){\tiny$s\geq2$}\put(18,10){\tiny$\cdots$}}}
\put(330,0){\put(0,0){\line(1,0){10}}\put(10,0){\line(0,1){50}}
\put(10,50){\line(-1,0){10}}\put(0,50){\line(0,-1){50}}
\put(1,-10){\tiny$a$}
\put(10,-20){\put(0,0){\line(1,0){10}}\put(10,0){\line(0,1){90}}
\put(10,90){\line(-1,0){10}}\put(0,90){\line(0,-1){90}}
\put(-15,-15){\tiny$s\geq a+2$}\put(18,36){\tiny$\cdots$}}}
\color{mygray}
\multiput(0.2,0.2)(0,.332){150}{\line(1,0){9.7}}
\multiput(50.2,0.2)(0,.332){150}{\line(1,0){9.7}}
\multiput(140.2,0.2)(0,.332){150}{\line(1,0){9.7}}
\multiput(240.2,0.2)(0,.332){150}{\line(1,0){9.7}}
\multiput(330.2,0.2)(0,.332){150}{\line(1,0){9.7}}
\end{picture}
\caption{{\it Decomposition} of a polyomino in $\CCP_a$.}\label{figdd}
\end{figure}
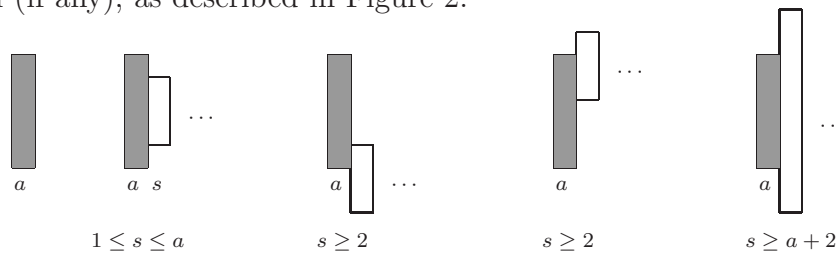

For the rest of this paper, we use $x,$ $y$, $p$, and $q$ to mark $h(.)$, $v(.)$, $B(.)$ and $U(.)$, respectively, while defining the corresponding generating functions. This section follows the methodology implemented in \cite{CakManYil} closely. Define $\CCP$ to be the set of all nonempty column-convex polyominoes. Similarly, set $\CCP_a$ to be the set of polyominoes in $\CCP$ with $a$ cells in their first columns. We define $F_a:=F_a(x,y,p,q)$ to be the generating function for the polyominoes in $\CCP_a$ according to $h(.)$, $v(.)$, $B(.)$, and $U(.)$; that is,
\beq
F_a(x,y,p,q) := \sum_{\nu \in \CCP_a} x^{h(\nu)}y^{v(\nu)}p^{B(\nu)}q^{U(\nu)}.
\feq
Next, decompose each $\nu\in \CCP_a$ by considering the size and the bottom position of its second column (if any), as described in Figure~\ref{figdd}.

By a careful analysis of these cases, for each $a\geq1$, we may write
\begin{align*}
F_a&=xy^a+\sum_{s=1}^{a-1}(a-1-s+p+q)xy^{a-s}F_s+pqxF_a\\
&+2\sum_{s=2}^ax\frac{y^{a+1-s}-y^a}{1-y}F_s+\sum_{s\geq a+1}x(2\frac{y-y^a}{1-y}+p+q)F_s
+\sum_{s\geq a+2}(s-1-a)xF_s.
\end{align*}
Now, define
\beq
F(u)=F(u;x,y,p,q)=\sum_{a\geq1}F_au^{a-1}.
\feq
By multiplying the last recurrence by $u^{a-1}$ and summing over all $a\geq1$, we obtain the following result.
\begin{lemma}\label{thc1}
The generating function $F(u)=F(u;x,y,p,q)$ satisfies
\begin{align*}
&K(u;x,y,p,q)F(u)\\
&=\frac{xy}{1-yu}+\frac{x}{1-u}\left(\frac{2yu}{1-yu}+p+q-\frac{1}{1-u}\right)F(1)+\frac{x}{1-u}\frac{\partial}{\partial u}F(u)\mid_{u=1},
\end{align*}
where $K(u;x,y,p,q)=1-\frac{(qyu(u-1)-yu^2-qu+2yu+q-1)(pyu(u-1)-yu^2-pu+2yu+p-1)x}{(1-u)^2(1-yu)^2}$.
\end{lemma}

Note that by setting  $q=p$ in Lemma \ref{thc1}, we get
\begin{align*}
&K(u;x,y,q,q)F(u)\\
&=\frac{xy}{1-yu}+\frac{x}{1-u}\left(\frac{2yu}{1-yu}+2q-\frac{1}{1-u}\right)F(1)+\frac{x}{1-u}\frac{\partial}{\partial u}F(u)\mid_{u=1},
\end{align*}
where
\beq
K(u;x,y,q,q)=1-\frac{(qyu(u-1)-yu^2-qu+2yu+q-1)^2x}{(1-u)^2(1-yu)^2}.
\feq
To solve this functional equation, we apply kernel method (see \cite{Ban}). To that end, setting $K(u;x,y,q,q)=0$ yields
$$\mp1+(q-1)\sqrt{x}+((\pm1-\sqrt{x}q)(y+1)+2\sqrt{x}y)u+(\mp y+(q-1)\sqrt{x}y)u^2=0.$$
Solving for $u$, the roots of this equation are
\begin{align*}
u_\pm&=\frac{1-\sqrt{1-\frac{4y((q\sqrt{x}\pm1)^2+\sqrt{x}(2q-1)(\mp2-\sqrt{x}))}
{((1\mp\sqrt{x}q)(y+1)\pm2\sqrt{x}y)^2}}}
{\frac{2y(1+\sqrt{x}(q-1))}{((1\mp\sqrt{x}q)(y+1)\pm2\sqrt{x}y)^2}}.
\end{align*}
By substituting $u=u_\pm$ into Lemma \ref{thc1}, solving for $F(1;x,y,q,q)$, and evaluating $\frac{\partial}{\partial u}F(u;x,y,q,q)\mid_{u=1}$, we obtain the main result of this section.
\begin{theorem}\label{thc2}
The generating function $F(1;x,y,q,q)$ is given by
\begin{align*}
&\frac{(u_+-1)(u_--1)y(y-1)}{u_+u_-y(y-2)+(u_++u_-)y-2y+1}\\
&=xy+xy^2+q^2x^2y+xy^3+( q^2+4q)x^2y^2+q^4x^3y+\cdots.
\end{align*}
In particular, $F(1;x,y,0,0)$ counts the number
of the column-convex Carlitz polyominoes according to $h(.)$ and $v(.)$, and is given by
\begin{align*}
&\frac{(v_+-1)(v_--1)y(y-1)}{v_+v_-y(y-2)+(v_++v_-)y-2y+1}\\
&=xy+xy^2+xy^3+xy^4+4x^2y^3+xy^5+12x^2y^4+x^3y^3+\cdots,
\end{align*}
where $v_\pm=u_\pm\mid_{q=0}$.
\end{theorem}
Note that the generating function for the number
of the column-convex Carlitz polyominoes according to the perimeter is simply given by
\beqn \label{rez1}
F(1;x^2,x^2,0,0)&=&
\frac{(1-x^2)(20x^6-45x^4+42x^2-21)}{4(8x^6-27x^4+36x^2-18)} \notag \\
&&+\frac{(1-x)(2x^3-9x^2+9)(x+1)^2}{4(8x^6-27x^4+36x^2-18)}\sqrt{4x^4+x^2+1-2x(1+2x^2)}\\
&&-\frac{(1+x)(2x^3+9x^2-9)(x-1)^2}{4(8x^6-27x^4+36x^2-18)}\sqrt{4x^4+x^2+1+2x(1+2x^2)}\notag \\
&&+\frac{3(x^2-1)^2}{4(8x^6-27x^4+36x^2-18)}\sqrt{(4x^4+x^2+1)^2-4x^2(1+2x^2)^2}. \notag
\feqn
We end this section by extracting the coefficient of $x^{n}$ in $F(1;x,x,0,0)$ and $\frac{d}{dq}F(1;x,x,q,q)\mid_{q=1}$, for large $n$, by conducting singularity analysis of \eqref{rez1}. The end result states that
\begin{corollary}\label{cor2.4}
Asymptotically, as $n$ gets large,
\begin{enumerate}
\item the number of column-convex Carlitz polyominoes with perimeter $2n$ is  $\frac{9\sqrt{2}(14+3\sqrt{3})}{2704\sqrt{\pi n^3}}4^n$.
\item  the total sum of $B+U$ over all column-convex polyominoes with perimeter $2n$ is
\beq
\frac{(1588-999\sqrt{2})\sqrt{5\sqrt{2}-7}
+6(51\sqrt{2}-28)\sqrt{99\sqrt{2}-140}}
{2209\sqrt{\pi n}}(\frac{3+2\sqrt{2}}{2})^n.
\feq
\end{enumerate}
\end{corollary}

\section{Enumeration of convex Carlitz polyominoes} \label{sec-b}

\begin{figure}[htp]
	\begin{picture}(55,45)
	\put(0,5){		\multiput(0,0)(0,4){9}{\put(0,0){\line(1,0){4}}\put(4,0){\line(0,1){4}}
			\put(4,4){\line(-1,0){4}}\put(0,4){\line(0,-1){4}}}		\multiput(4,-8)(0,4){9}{\put(0,0){\line(1,0){4}}\put(4,0){\line(0,1){4}}
			\put(4,4){\line(-1,0){4}}\put(0,4){\line(0,-1){4}}} \multiput(8,-4)(0,4){6}{\put(0,0){\line(1,0){4}}\put(4,0){\line(0,1){4}}
			\put(4,4){\line(-1,0){4}}\put(0,4){\line(0,-1){4}}} \multiput(12,-4)(0,4){6}{\put(0,0){\line(1,0){4}}\put(4,0){\line(0,1){4}}
			\put(4,4){\line(-1,0){4}}\put(0,4){\line(0,-1){4}}} \multiput(16,4)(0,4){3}{\put(0,0){\line(1,0){4}}\put(4,0){\line(0,1){4}}
			\put(4,4){\line(-1,0){4}}\put(0,4){\line(0,-1){4}}} }
	\put(40,5){
		\multiput(0,0)(0,4){9}{\put(0,0){\line(1,0){4}}\put(4,0){\line(0,1){4}}
			\put(4,4){\line(-1,0){4}}\put(0,4){\line(0,-1){4}}}		\multiput(4,0)(0,4){7}{\put(0,0){\line(1,0){4}}\put(4,0){\line(0,1){4}}
			\put(4,4){\line(-1,0){4}}\put(0,4){\line(0,-1){4}}}
		\multiput(8,4)(0,4){6}{\put(0,0){\line(1,0){4}}\put(4,0){\line(0,1){4}}
			\put(4,4){\line(-1,0){4}}\put(0,4){\line(0,-1){4}}} \multiput(12,8)(0,4){4}{\put(0,0){\line(1,0){4}}\put(4,0){\line(0,1){4}}
			\put(4,4){\line(-1,0){4}}\put(0,4){\line(0,-1){4}}} \multiput(16,12)(0,4){2}{\put(0,0){\line(1,0){4}}\put(4,0){\line(0,1){4}}
			\put(4,4){\line(-1,0){4}}\put(0,4){\line(0,-1){4}}} }
	\end{picture}
	\caption{An example of polyominoes in (left) $\CP^t$ (right) $\CP^{bt}$  \label{fig10}}
\end{figure}
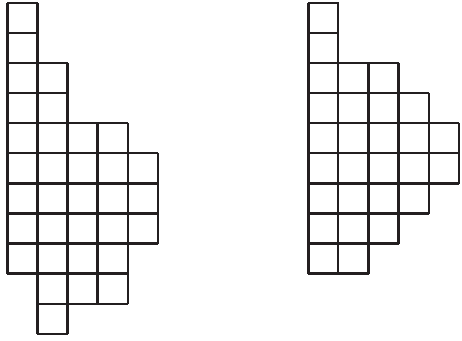

This section follows the approach developed in \cite{MR} closely. We define $\CP^t$ to be the set of all nonempty convex polyominoes $\nu,$ where for each column $j$ of $\nu$, $u_s(\nu)\leq u_j(\nu)$ for all $s\geq j+1$. Similarly, $\CP^b$ is the set of all nonempty convex polyominoes $\nu$ such that for all columns $j,$ and for all $s\geq j+1,$ $b_s(\nu)\geq b_j(\nu)$. Set $\CP^{bt}:=\CP^t\cap\CP^b$.  By a clear upside down symmetry, there is a bijection between the set $\CP^t$ and the set $\CP^b$. See Figure \ref{fig10}, for an example of polyominoes in $\CP^t$ and $\CP^{bt}$. In addition, we let $\CP_a^t,$ $\CP_a^b,$ $\CP_a^{bt},$ and $\CP_a$ denote the set of polyominoes in $\CP^t,$ $\CP^b,$ $\CP^{bt},$ and $\CP$ with $a$ cells in their first columns. Our enumeration is done as follows. The first step is to count the number of polyominoes in $\CP^{bt}$ with respect to the statistics of interest $h(.)$, $v(.)$, $U(.)$, and $B(.)$. Then, for the second step, we extend the result from $\CP^{bt}$ to $\CP^{t}$ and $\CP^{b}.$ This can be done since the corresponding generating functions can be written recursively in terms of the corresponding generating functions in $\CP^{bt}$. When this is done, the last step is to obtain the results for $\CP$ by lifting up the result obtained for $\CP^t$ and $\CP^b$.

Let $G_a:=G_a(x,y,p,q)$ be the generating functions of $h(.),$ $v(.),$ $B(.)$, and $U(.)$ over $\CP_a$; that is,
\beq
G_a:=G_a(x,y,p,q) = \sum_{\nu\in \CP_a} x^{h(\nu)}y^{v(\nu)}p^{B(\nu)}q^{U(\nu)}.
\feq
We define $G_a^t$, $G_a^b$, $G_a^{bt}$ similarly for the set of polyominoes $\CP_a^b$, $\CP_a^t$, $\CP_a^{bt}$, respectively. Moreover, we define the generating function $G(u)=G(u;x,y,p,q)=\sum_{a\geq1}G_au^{a-1}$. Additionally, we define $G^b(u),G^t(u)$ and $G^{bt}(u)$.

Consider a polyomino $\nu$ in $\CP_a^{bt}$ with the second column of size $s$ (if any). If second column does not exist then the contribution to $G_a^{bt}$ is exactly $xy^a$. If second column exists then $u_2(\nu)\leq u_1(\nu)$ and $b_2(\nu)\geq b_1(\nu)$, so there are exactly $a+1-s$ ways of gluing the second column to the first column. By considering all four cases
\begin{enumerate}
	\item $u_2(\nu)<u_1(\nu),b_2(\nu)>b_1(\nu)$
	\item $u_2(\nu)<u_1(\nu),b_2(\nu)=b_1(\nu)$,
	\item $u_2(\nu)=u_1(\nu),b_2(\nu)>b_1(\nu)$, and
	\item $u_2(\nu)=u_1(\nu),b_2(\nu)=b_1(\nu)$,
\end{enumerate}
we obtain
\beq
G_a^{bt}=xy^a+\sum_{s=1}^{a-1}(a-1-s+p+q)xy^{a-s}G_s^{bt}+pqxG_a^{bt}.
\feq
By multiplying by $u^{a-1}$ and summing over $a\geq1$, we obtain
$$G^{bt}(u)=\frac{xy}{1-yu}+\frac{xyu(-yu(p+q-1)+p+q)}{(1-yu)^2}
G^{bt}(u)+pqxG^{bt}(u).$$
Hence,
\begin{align}\label{eq1Gbu}
G^{bt}(u;x,y,p,q)=\frac{\frac{xy}{1-yu}}{1-pqx-\frac{xyu(-yu(p+q-1)+p+q)}{(1-yu)^2}}.
\end{align}

\begin{figure}[htp]
\begin{picture}(80,26)(-20,-10)
\setlength{\unitlength}{.3mm}
\put(-100,0){\put(0,0){\line(1,0){10}}\put(10,0){\line(0,1){50}}
\put(10,50){\line(-1,0){10}}\put(0,50){\line(0,-1){50}}
\put(1,-10){\tiny$a$}}
\put(-50,0){\put(0,0){\line(1,0){10}}\put(10,0){\line(0,1){50}}
\put(10,50){\line(-1,0){10}}\put(0,50){\line(0,-1){50}}
\put(1,-10){\tiny$a$}
\put(10,10){\put(0,0){\line(1,0){10}}\put(10,0){\line(0,1){30}}
\put(10,30){\line(-1,0){10}}\put(0,30){\line(0,-1){30}}
\put(2,-20){\tiny$s$}\put(18,10){\tiny$\cdots$}}
\put(-25,-25){\tiny$1\leq s\leq a-2$}}
\put(0,0){\put(0,0){\line(1,0){10}}\put(10,0){\line(0,1){50}}
\put(10,50){\line(-1,0){10}}\put(0,50){\line(0,-1){50}}
\put(1,-10){\tiny$a$}
\put(10,10){\put(0,0){\line(1,0){10}}\put(10,0){\line(0,1){40}}
\put(10,40){\line(-1,0){10}}\put(0,40){\line(0,-1){40}}
\put(2,-20){\tiny$s$}\put(18,10){\tiny$\cdots$}}
\put(-15,-35){\tiny$1\leq s\leq a-1$}}
\put(50,0){\put(0,0){\line(1,0){10}}\put(10,0){\line(0,1){50}}
\put(10,50){\line(-1,0){10}}\put(0,50){\line(0,-1){50}}
\put(1,-10){\tiny$a$}
\put(10,0){\put(0,0){\line(1,0){10}}\put(10,0){\line(0,1){40}}
\put(10,40){\line(-1,0){10}}\put(0,40){\line(0,-1){40}}
\put(2,-10){\tiny$s$}\put(18,10){\tiny$\cdots$}}
\put(-25,-25){\tiny$1\leq s\leq a-1$}}
\put(100,0){\put(0,0){\line(1,0){10}}\put(10,0){\line(0,1){50}}
\put(10,50){\line(-1,0){10}}\put(0,50){\line(0,-1){50}}
\put(1,-10){\tiny$a$}
\put(10,0){\put(0,0){\line(1,0){10}}\put(10,0){\line(0,1){50}}
\put(10,50){\line(-1,0){10}}\put(0,50){\line(0,-1){50}}
\put(2,-10){\tiny$a$}\put(18,10){\tiny$\cdots$}}}
\put(150,0){\put(0,0){\line(1,0){10}}\put(10,0){\line(0,1){50}}
\put(10,50){\line(-1,0){10}}\put(0,50){\line(0,-1){50}}
\put(1,-10){\tiny$a$}
\put(10,-20){\put(0,0){\line(1,0){10}}\put(10,0){\line(0,1){30}}
\put(10,30){\line(-1,0){10}}\put(0,30){\line(0,-1){30}}
\put(-15,-15){\tiny$s\geq2$}\put(18,10){\tiny$\cdots$}}}
\put(200,0){\put(0,0){\line(1,0){10}}\put(10,0){\line(0,1){50}}
\put(10,50){\line(-1,0){10}}\put(0,50){\line(0,-1){50}}
\put(1,-10){\tiny$a$}
\put(10,-20){\put(0,0){\line(1,0){10}}\put(10,0){\line(0,1){70}}
\put(10,70){\line(-1,0){10}}\put(0,70){\line(0,-1){70}}
\put(-15,-15){\tiny$s\geq2$}\put(18,10){\tiny$\cdots$}}}
\color{mygray}
\multiput(-99.8,0.2)(0,.332){150}{\line(1,0){9.7}}
\multiput(-49.8,0.2)(0,.332){150}{\line(1,0){9.7}}
\multiput(0.2,0.2)(0,.332){150}{\line(1,0){9.7}}
\multiput(50.2,0.2)(0,.332){150}{\line(1,0){9.7}}
\multiput(100.2,0.2)(0,.332){150}{\line(1,0){9.7}}
\multiput(150.2,0.2)(0,.332){150}{\line(1,0){9.7}}
\multiput(200.2,0.2)(0,.332){150}{\line(1,0){9.7}}
\end{picture}
\caption{Decomposition of a polyomino in $\CP_a^t$.}\label{figuu}
\end{figure}
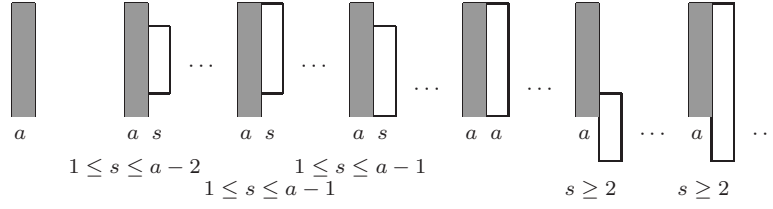

Now, let us find a formula for $G^t(u)$.
Any nonempty polyomino $\nu\in \CP^t_a$ can be decomposed as described in Figure~\ref{figuu}.
Therefore, the generating function $G^t_a(u)$ satisfies
\begin{align*}
G^t_a&=xy^a+\sum_{s=1}^{a-1}(a-1-s)xy^{a-s}G^{bt}_s
+q\sum_{s=1}^{a-1}xy^{a-s}G_s^{bt}+p\sum_{s=1}^{a-1}xy^{a-s}G_s^t\\
&+pqxG_a^t+\sum_{s=2}^ax\frac{y^{a+1-s}-y^a}{1-y}G_s^t
+\sum_{s\geq a+1}x\left(\frac{y-y^a}{1-y}+q\right)G_s^t.
\end{align*}
By multiplying this recurrence by $u^{a-1}$ and summing over $a\geq1$, we obtain the equation
\begin{align}
&\left(1-pqx+\frac{qx}{1-u}-\frac{pxyu}{1-yu}+\frac{xyu}{(1-u)(1-yu)}\right)G^t(u)\nonumber\\
&=\frac{xy}{1-yu}+\left(\frac{xy^2u^2}{(1-yu)^2}
+\frac{qxyu}{1-yu}\right)G^{bt}(u)+\frac{xyu+qx(1-yu)}{(1-u)(1-yu)}G^t(1),\label{eq1Gu}
\end{align}
which will be solved by using the kernel method (see the beginning of Section \ref{sec-a}). First, by finding the roots of the equation
$$1-pqx+\frac{qx}{1-u'}-\frac{pxyu'}{1-yu'}+\frac{xyu'}{(1-u')(1-yu')}=0$$
for variable $u$, we get
\beq
u'=\frac{1-\sqrt{1-\frac{4y(1+p(1-q)x)(1+q(1-p)x)}{(1+y-pqx-xy(1-p)(1-q))^2}}}
{\frac{2y(1+p(1-q)x)}{1+y-pqx-xy(1-p)(1-q)}}.
\feq
Then, setting $u=u'$ in last functional equation yields
\begin{align*}
G^t(1)=\frac{y(u'-1)}{yu'+q(1-yu')}+\frac{yu'(u'-1)}{1-yu'}G^{bt}(u').
\end{align*}
Thus, by \eqref{eq1Gbu}, we can state the following result.
\begin{lemma}\label{lem1u}
The generating function $G^t(u;x,y,p,q)$ is given by \eqref{eq1Gu}, where the generating function $G^t(1;x,y,p,q)$ is given by
\begin{align}\label{eq1Gu1}
G^t(1)&=\frac{y(u'-1)}{yu'+q(1-yu')}+\frac{xy^2u'(u'-1)}
{(1-pqx)(1-yu')^2-xyu'(p+q)(1-yu')-xy^2u'^2}.
\end{align}
\end{lemma}

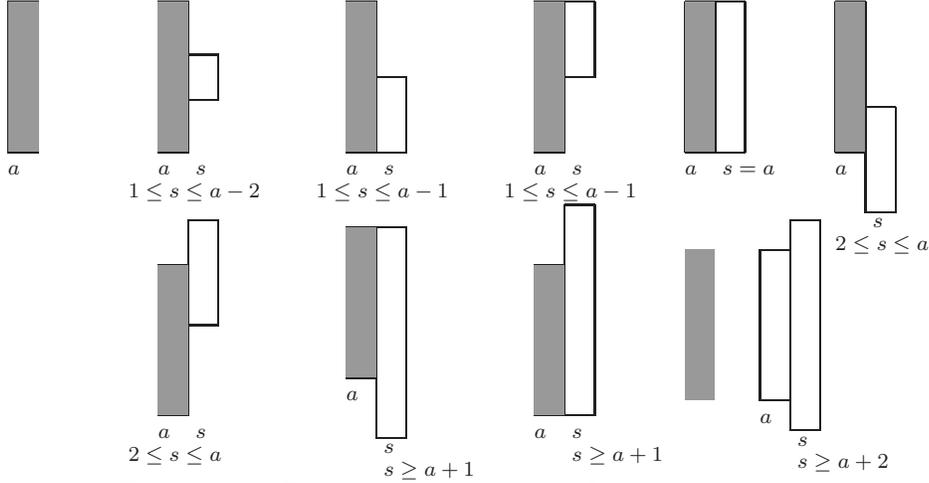
\begin{figure}[htp]
\begin{picture}(120,60)
	\setlength{\unitlength}{0.1cm} \put(0,-7){
\put(0,30){
		\put(0,15){ \put(0,0){\line(1,0){4}}\put(4,0){\line(0,1){20}}\put(4,20){\line(-1,0){4}}\put(0,20){\line(0,-1){20}}\put(0,-3){\tiny$a$}}
		\put(20,15){ \put(0,0){\line(1,0){4}}\put(4,0){\line(0,1){20}}\put(4,20){\line(-1,0){4}}\put(0,20){\line(0,-1){20}}		\put(4,7){\put(0,0){\line(1,0){4}}\put(4,0){\line(0,1){6}}\put(4,6){\line(-1,0){4}}\put(0,6){\line(0,-1){6}}}
\put(0,-3){\tiny$a$}\put(5,-3){\tiny$s$}\put(-4,-6){\tiny$1\leq s\leq a-2$}}
		\put(45,15){ \put(0,0){\line(1,0){4}}\put(4,0){\line(0,1){20}}\put(4,20){\line(-1,0){4}}\put(0,20){\line(0,-1){20}} \put(4,0){\put(0,0){\line(1,0){4}}\put(4,0){\line(0,1){10}}\put(4,10){\line(-1,0){4}}\put(0,10){\line(0,-1){10}}}
\put(0,-3){\tiny$a$}\put(5,-3){\tiny$s$}\put(-4,-6){\tiny$1\leq s\leq a-1$}}
		\put(70,15){ \put(0,0){\line(1,0){4}}\put(4,0){\line(0,1){20}}\put(4,20){\line(-1,0){4}}\put(0,20){\line(0,-1){20}}			\put(4,10){\put(0,0){\line(1,0){4}}\put(4,0){\line(0,1){10}}\put(4,10){\line(-1,0){4}}\put(0,10){\line(0,-1){10}}}
\put(0,-3){\tiny$a$}\put(5,-3){\tiny$s$}\put(-4,-6){\tiny$1\leq s\leq a-1$}}
		\put(90,15){ \put(0,0){\line(1,0){4}}\put(4,0){\line(0,1){20}}\put(4,20){\line(-1,0){4}}\put(0,20){\line(0,-1){20}} \put(4,0){\put(0,0){\line(1,0){4}}\put(4,0){\line(0,1){20}}\put(4,20){\line(-1,0){4}}\put(0,20){\line(0,-1){20}}}
\put(0,-3){\tiny$a$}\put(5,-3){\tiny$s=a$}}
		\put(110,15){ \put(0,0){\line(1,0){4}}\put(4,0){\line(0,1){20}}\put(4,20){\line(-1,0){4}}\put(0,20){\line(0,-1){20}}			\put(4,-8){\put(0,0){\line(1,0){4}}\put(4,0){\line(0,1){14}}\put(4,14){\line(-1,0){4}}\put(0,14){\line(0,-1){14}}}
\put(0,-3){\tiny$a$}\put(5,-10){\tiny$s$}\put(0,-13){\tiny$2\leq s\leq a$}}
}
\put(-110,0){
		\put(130,10){	\put(0,0){\line(1,0){4}}\put(4,0){\line(0,1){20}}\put(4,20){\line(-1,0){4}}\put(0,20){\line(0,-1){20}} \put(4,12){\put(0,0){\line(1,0){4}}\put(4,0){\line(0,1){14}}\put(4,14){\line(-1,0){4}}\put(0,14){\line(0,-1){14}}}
\put(0,-3){\tiny$a$}\put(5,-3){\tiny$s$}\put(-4,-6){\tiny$2\leq s\leq a$}}
		\put(155,15){ \put(0,0){\line(1,0){4}}\put(4,0){\line(0,1){20}}\put(4,20){\line(-1,0){4}}\put(0,20){\line(0,-1){20}}			\put(4,-8){\put(0,0){\line(1,0){4}}\put(4,0){\line(0,1){28}}\put(4,28){\line(-1,0){4}}\put(0,28){\line(0,-1){28}}}
\put(0,-3){\tiny$a$}\put(5,-10){\tiny$s$}\put(5,-13){\tiny$s\geq a+1$}}
		\put(180,10){ \put(0,0){\line(1,0){4}}\put(4,0){\line(0,1){20}}\put(4,20){\line(-1,0){4}}\put(0,20){\line(0,-1){20}} \put(4,0){\put(0,0){\line(1,0){4}}\put(4,0){\line(0,1){28}}\put(4,28){\line(-1,0){4}}\put(0,28){\line(0,-1){28}}}
\put(0,-3){\tiny$a$}\put(5,-3){\tiny$s$}\put(5,-6){\tiny$s\geq a+1$}}
		\put(210,12){ \put(0,0){\line(1,0){4}}\put(4,0){\line(0,1){20}}\put(4,20){\line(-1,0){4}}\put(0,20){\line(0,-1){20}}			\put(4,-4){\put(0,0){\line(1,0){4}}\put(4,0){\line(0,1){28}}\put(4,28){\line(-1,0){4}}\put(0,28){\line(0,-1){28}}}
\put(0,-3){\tiny$a$}\put(5,-6){\tiny$s$}\put(5,-9){\tiny$s\geq a+2$}} }}
\color{mygray}
\multiput(0.1,38.2)(0,.1){198}{\line(1,0){3.9}}
\multiput(20,38.2)(0,.1){198}{\line(1,0){3.9}}
\multiput(45,38.2)(0,.1){198}{\line(1,0){3.9}}
\multiput(70,38.2)(0,.1){198}{\line(1,0){3.9}}
\multiput(90,38.2)(0,.1){198}{\line(1,0){3.9}}
\multiput(110,38.2)(0,.1){198}{\line(1,0){3.9}}
\multiput(20,3.2)(0,.1){198}{\line(1,0){3.9}}
\multiput(45,8.2)(0,.1){198}{\line(1,0){3.9}}
\multiput(70,3.2)(0,.1){198}{\line(1,0){3.9}}
\multiput(90,5.2)(0,.1){198}{\line(1,0){3.9}}
	\end{picture}
	\caption{Decomposition of a polyomino in $\CP_a$}\label{figcpcp}
\end{figure}

Next, we find the generating function $G(1;x,y,p,q)$ based on $G^t$ and $G^b$. To that goal, we decompose each convex polyomino as described in Figure \ref{figcpcp}.
Considering these cases, we can rewrite $G_a(x,y,t)$ as
\begin{align*}
G_a&=xy^a+\sum_{s=1}^{a-1}(a-1-s)xy^{a-s}G_b^{bt}
+p\sum_{s=1}^{a-1}xy^{a-s}G_b^{t}+q\sum_{s=1}^{a-1}xy^{a-s}G_s^{b}\\
&+xG_a+\sum_{s=2}^ax(y^{a-1}+\cdots+y^{a+1-s})G_s^t+\sum_{s=2}^ax(y^{a-1}+\cdots+y^{a+1-s})G_s^b\\
&+\sum_{s\geq a+1}x(y^{a-1}+\cdots+y)G_s^t+\sum_{s\geq a+1}x(y^{a-1}+\cdots+y)G_s^b+(p+q)\sum_{s\geq a+1}xG_s\\
&+\sum_{s\geq a+1}(s-1-a)xG_s.
\end{align*}
By multiplying by $u^{a-1}$ and summing over $a\geq1$, we obtain
\beqn \label{res2}
&\left(1-pqx+\frac{(p+q)x}{1-u}-\frac{x}{(1-u)^2}\right)G(u)\notag \\
&\qquad=\frac{xy}{1-yu}+\frac{xy^2u^2}{(1-yu)^2}G^{bt}(u)+\frac{xyu}{1-yu}(pG^{t}(u)+qG^{b}(u))\notag \\
&\qquad+\frac{xyu}{(1-u)(1-yu)}(G^t(1)+G^b(1)-G^t(u)-G^b(u))\notag \\
&\qquad+\frac{(p+q)x}{1-u}G(1)-\frac{x}{(1-u)^2}G(1)+\frac{x}{1-u}\frac{\partial}{\partial u}G(u)\mid_{u=1}.
\feqn
Recall $G^b(u;x,y,p,q)=G^t(u;x,y,q,p)$ by a symmetry argument. Also, note that the roots of the kernel
\beq
K(u)=1-pqx+\frac{(p+q)x}{1-u}-\frac{x}{(1-u)^2}
\feq
are
\beq
u_\pm=1+\frac{(p+q)x\pm\sqrt{(p-q)^2x^2+4x}}{2(1-pqx)}.
\feq
By inserting $u=u_+$ and $u=u_-$, one at the time, into the functional equation \eqref{res2}, we obtain a system of two equations with variables $G(1)$ and $\frac{\partial}{\partial u}G(u)\mid_{u=1}$. Solving this system yields

\begin{theorem}\label{thMG1}
The generating function $G(1;x,y,p,q)$ is given by
\begin{align*}
&\frac{y(y-1)(1-u_+)(1-u_-)}{(1-yu_+)(1-yu_-)}\\
&+\frac{y^2(1-u_+)(1-u_-)}{u_+-u_-}
\left(\frac{u_+^2(1-u_+)G^{bt}(u_+;x,y,p,q)}{(1-yu_+)^2}
-\frac{u_-^2(1-u_-)G^{bt}(u_-;x,y,p,q)}{(1-yu_-)^2}\right)\\
&-\frac{yu_+(1-u_+)(1-u_-)}{(u_+-u_-)(1-yu_+)}\left((1-p(1-u_+))G^t(u_+;x,y,p,q)
+(1-q(1-u_+))G^t(u_+;x,y,q,p)\right)\\
&+\frac{yu_-(1-u_+)(1-u_-)}{(u_+-u_-)(1-yu_-)}\left((1-p(1-u_-))G^t(u_-;x,y,p,q)
+(1-q(1-u_-))G^t(u_-;x,y,q,p)\right)\\
&+\frac{y(1-u_+)(1-u_-)}{(1-yu_+)(1-yu_-)}(G^t(1;x,y,p,q)+G^b(1;x,y,q,p)),
\end{align*}
where $G^{bt}(u;x,y,p,q)$, $G^t(u;x,y,p,q)$ and $G^t(1;x,y,p,q)$  are given by \eqref{eq1Gbu},  \eqref{eq1Gu} and \eqref{eq1Gu1}, respectively.
\end{theorem}
As a corollary, the generating function $G(1;x,x,q,q)$ is given by
\begin{align*}
&\frac{x^2A}{(q^2x^2-q^2x-2qx^2+x^2-3x+1)^2(q^2x^2-q^2x-2qx^2+x^2+x+1)^2}\\
&-\frac{x^4(q^3x^2-q^3x-3q^2x^2+q^2x+3qx^2-qx-x^2+q+x+1)^2}{
(q^2x^2-q^2x-2qx^2+x^2-3x+1)^{3/2}(q^2x^2-q^2x-2qx^2+x^2+x+1)^{3/2}}\\
&=x^2+(q^2+1)x^3+(q^4+q^2+4*q+1)x^4+(q^2+1)(q^4+8q+5)x^5+\cdots,
\end{align*}
where
\begin{align*}
A&=1-3(q^2+1)x+q(3q^3+10q-2)x^2-(q^2+1)(q^4+11q^2-6q-4)x^3\\
&+(q-1)(6q^5+17q^3-7q^2-5q-3)x^4-(q^2+1)(q^4+11q^2-6q-4)(q-1)^2x^5\\
&+q(3q^3+10q-2)(q-1)^4x^6-3(q^2+1)(q-1)^6x^7+(q-1)^8x^8.
\end{align*}

By considering $G(1;x,x,0,0)$, we show that the generating function for the number of convex Carlitz  polyominoes according to their perimeters is given by
$$\frac{x^4(x^{16}-3x^{14}+4x^{10}+3x^8+4x^6-3x^2+1)}
{(x^4-3x^2+1)^2(x^4+x^2+1)^2}-\frac{x^8(x^4-x^2-1)^2}{(x^4-3x^2+1)^{3/2}(x^4+x^2+1)^{3/2}}.$$
Next, by singularity analysis of $G(1;x,x,0,0)$, we get
\begin{corollary}\label{cor3.3}
As $n$ grows to infinity, the number of convex Carlitz polyominoes with the perimeter $2n$ is asymptotic to $\frac{n+1}{10}\left(\frac{3+\sqrt{5}}{2}\right)^{n-2}$.
\end{corollary}
Note that by Theorem \ref{thMG1},
\begin{align*}
\frac{d}{dq}G(1;x,x,q,q)\mid_{q=1}
&=\frac{2x^3(16x^4-12x^3+10x^2-5x+1)}{(1-4x)^3}
+\frac{4x^4(4x^2-1)}{(1-4x)^{5/2}}.
\end{align*}
We conclude this section by stating that singularity analysis of $\frac{d}{dq}G(1;x,x,q,q)\mid_{q=1}$ also yields
\begin{corollary}\label{cor3.4}
The total sum of $B+U$ over all convex polyominoes with perimeter $2n$ is asymptotic to $n^24^{n-4}$.
\end{corollary}

\end{document}